\documentclass{article}                     

\usepackage{graphicx}
\usepackage{algorithm,algorithmicx,algpseudocode}
\newtheorem{defi}{Definition}
\usepackage{lineno,hyperref}
\usepackage{amsmath}

\begin{document}

\title{$L_1$-regularization for multi-period portfolio selection}


\renewcommand{\thefootnote}{\fnsymbol{footnote}}

\author{S. Corsaro\footnotemark[1],
        V. De Simone\footnotemark[2],
        Z. Marino\footnotemark[1],
        F. Perla\footnotemark[1]}

\date{}
\maketitle

\footnotetext[1]{Department of Management and Quantitative Studies, University of Naples ``Parthenope'',
              email: stefania.corsaro,zelda.marino,francesca.perla@uniparthenope.it}
\footnotetext[2]{Department of Mathematics and Physics, University of Campania ``Luigi Vanvitelli'',
email: valentina.desimone@unicampania.it}
\renewcommand{\thefootnote}{\arabic{footnote}}

\begin{abstract}
In this work we present a model for the solution of the multi-period portfolio selection problem.
The model is based on a time consistent dynamic risk measure.
We apply $l_1$-regularization to stabilize
the solution process and to obtain sparse solutions, which allow one to reduce holding costs.
The core problem is a nonsmooth optimization one, with
equality constraints.
We present an iterative procedure based on a modified Bregman iteration,
that adaptively sets the value of the regularization parameter in order
to produce solutions with desired financial properties.
We validate the approach showing results of tests performed on real data.\\
\textbf{Keywords}: {Portfolio optimization; Time consistency; $l_1$-norm; Constrained optimization}
\end{abstract}

\section{Introduction}
\label{intro}
In this work we focus on dynamic portfolio selection problem.
This problem arises in medium and long-term investments,
in which one allows decisions to change over time by the end of the investment,
taking into account the time evolution of available information.
We consider dynamic decision problems formulated in a discrete multi-stage setting,
with underlying time evolving continuously.
This formulation is usually referred to as multi-period portfolio selection.
In \cite{Li}, one of the pioneering works in this framework,
authors present an extension of Markowitz mean-variance model
to multi-period investments. However, their model does not satisfy
the \emph{time consistency} property, which has been recognized to
play a relevant role. This concept deals with the consistency over time
of decisions taken with the support of evolving available information.
Different definitions of time consistency can be found in literature,
either related to risk measures or investment policies \cite{Consigli}.
We consider the first case, in which one focuses on the properties of the multi-period
risk measure employed for modelling the investment problem;
this is time consistent if, according to it, the time evolving filtration
related to the evaluation of a stochastic process does not modify
decisions taken using values computed previously.
In other words, if we today establish that two investments
have the same level of risk,
then the same level of riskiness should have been estimated
for them yesterday.\\
In \cite{Chen} authors introduce a a time consistent mean-variance
model, which extends the classical Markowitz approach to the multi-period case.
The model is based on a separable expected conditional mapping, obtained
by summing single-period terms; it is shown that this ensures that
the underlying risk measure is time consistent.\\
We then adopt a separable formulation,
in which single-period terms are defined taking the variance
as risk measure; more precisely, we consider single-period minimum variance terms
and fix a final target expected wealth. Thus, we consider medium and
long-term investments in which investors are
not interested with positions at intermediate
periods.\\
A common choice to estimate Markowitz model parameters is to use historical data;
correlation between assets returns can lead to
ill-conditioned covariance matrices.
Different regularization techniques have been suggested
with the aim of improving the problem conditioning. Among these,
we consider $l_1$-regularization techniques to stabilize
the solution process.
The $l_1$-regularization for Markowitz model was
introduced in \cite{brodie}, where a
$l_1$-penalty term is added to promote sparsity in the solution.
Since solutions establish the amount of capital to be invested in each
available security,
sparsity means that money are invested in a few securities,
the so-called \emph{active positions}.
This allows investor to reduce the number of positions to be monitored and
the holding costs.
In this paper we propose a multi-period model expressed as a nonsmooth optimization problem with
equality constraints. One method to solve it is Bregman iteration \cite{Bregman},
that was recently introduced with success in many fields formulated in
term of $l_1$-minimization (image analysis \cite{Antonelli,GoldsteinOsher,compressed},
matrix rank \cite{rank}, portfolio selection \cite{Noi,maf18,Boyd}.
Following the idea presented in \cite{Noi} for the single-period case,
we present a modified version of Bregman iteration.
The method is modified to adaptively select the regularization parameter
that realizes a good trade-off between the fidelity to data and the financial properties
required to solution, given in terms of sparsity and number
of short positions.\\
Bregman iteration converts the constraint optimization problem in a short sequence of unconstrained ones;
The presence of the $l_1$-term makes the solution of the
involved optimization sub-problem not trivial, thus we
apply ad hoc methods to deal with non-smoothness \cite{fista}.\\
In section \ref{sec:2} we describe the regularized portfolio
selection model; in section \ref{sec:3} we describe Bregman iteration method.
In section \ref{sec:4} we show some numerical experiments performed on real data.

\section{Regularized portfolio selection model}
\label{sec:2}
In this section we introduce an $l_1$-regularized model for
multi-period portfolio selection.
Let $m$ be the number of investment periods.
Decisions are assumed at the rebalancing dates $j, \; j=1, \ldots, m$; decision taken
at time $j$ is kept in the period $[j,j+1)$.
Let $1\leq t\leq m$ be the evaluation date. A conditional
risk mapping $\rho_{t,m}$ computes the risk of a multi-period loss process
at time $t$, that is, its value is the risk estimate at time $t$ of future losses.
Dynamic time consistency ensures that investors’ preferences
remain consistent over time.
In probabilistic terms, we have that the information is described by a
filtration that evolves according to the dynamics of
a certain underlying random process whose risk is to be estimated.
The dynamics is supposed
not to affect previous estimates under a dynamic time consistent risk measure.
To simplify, given two loss processes
that produce the same losses up to the evaluation date $t$,
if they are estimated equally risky at time $t$,
then the same is observed previously
under a dynamic consistent risk measure.
Time consistency has been recognized to be a relevant property in
dynamic asset allocation. Indeed, much effort has been addressed to
the analysis and proposal of time consistent risk measures;
we address reader to \cite{Consigli} and references therein for an overview
on the existing literature.
In \cite{Chen} a time consistent dynamic mean-variance
model is presented. It is based on a separable expected conditional mapping,
which authors show to be time consistent.
We here report the definition given in \cite{Chen}:
\begin{defi}
A conditional risk measure $\rho_{t,m}$ on the time horizon T=$[t,m]$
is separable if it can be
expressed in the following way:
\begin{equation}\label{eq:sep}
\rho_{t,m}(Y_{t,m}) = \sum_{s=t+1}^m E[\rho_s(Y_s|F_{s-1})|F_t]
\end{equation}
where $Y_{t,m}$ is a random vector of losses estimated at the dates
$\{t,t+1,\ldots,m\}$, $F_t$ the
filtration at time $t$,
$Y_t$ is the random loss adapted to $F_t$ and $E$ is the expectation operator.
\end{defi}
Relation \eqref{eq:sep}
states that the risk measure is decomposed into a sum of terms
in such a way that each one provides a risk estimate
in one period of the investment, using information available
at the beginning of the period.
In this paper we consider a separable conditional risk measure,
taking the variance as single-period risk measure.\\
Let us denote
the number of traded assets by $n$ and
with $\mathbf{w}_j$ the portfolio of holdings at the beginning
of period $j$.
Thus, for instance, $(\mathbf{w}_j)_i$ is the amount invested in the asset
$i$ at the $j$-th rebalancing date. The optimal portfolio is then
defined by the vector $\mathbf{w} = (\mathbf{w}_1, \mathbf{w}_2, \ldots, \mathbf{w}_m)\in \Re^N$, where
$N = m\cdot n$ is the problem dimension.
Finally, the vector $\mathbf{r}_j$ and
$C_j \in \Re^{n\times n}$ contain respectively the expected return vector and the covariance
matrix estimated at time $j$;
covariance matrices are assumed to be positive definite.
We aim at minimizing the risk of the strategy, estimated by
\begin{equation}\label{eq:fo}
\Phi(w) = \sum_{j=1}^{m} \mathbf{w}_j^T C_j\mathbf{w}_j.
\end{equation}
The objective function \eqref{eq:fo} typically leads to
ill-conditioned problems, because of assets correlation.
At this purpose, we apply $l_1$-regularization to
stabilize the solution process; an $l_1$-penalty term is added to \eqref{eq:fo}.
This technique was applied in the single-period case, that is, in the
classical Markowitz approach, in \cite{Noi}. In that paper,
authors show that $l_1$-regularization provides sparse solutions;
since solutions establish the amount of capital to be invested in each
available security, sparsity means that money are invested in a few securities.
This allows investor to reduce both the number of positions to
be monitored and the holding costs.

Thus, we propose the following constrained
optimization problem for multi-period portfolio selection:
\begin{equation}\label{eq:discreto}
\begin{array}{l}
    \min_\mathbf{w} \sum_{j=1}^{m} \Bigl[ \mathbf{w}_j^T C_j\mathbf{w}_j + \tau \|\mathbf{w}_j\|_1 \Bigr] \\
    \mathrm{s.t.} \\
    \mathbf{w}_1^T \mathbf{1}_n = \xi_{\mathrm{init}} \\
    \mathbf{w}_j^T \mathbf{1}_n = (\mathbf{1}_n+\mathbf{r}_{j-1})^T \mathbf{w}_{j-1}, \quad j = 2,\ldots,m \\
    \mathbf{w}_m^T \mathbf{1}_n = \xi_{\mathrm{term}} \\
\end{array}
\end{equation}
where $\tau>0$ is the regularization parameter,
$\xi_{\mathrm{init}}$ is the initial wealth, $\xi_{\mathrm{term}}$
is the target expected wealth resulting from the overall investment,
$\mathbf{1}_n$ is vector
of ones of length $n$.
The first constraint is the budget constraint.
We assume that the strategy is self-financed: this is stated in constraints from
$2$ to $m$, where it is established that
at the end of each period the wealth
is given by the revaluation of the previous one.
The $(m+1)$-th constraint defines the expected final wealth.

\section{Bregman iteration for portfolio selection}
\label{sec:3}
In this section we discuss the solution of the $l_1$-regularized
constrained optimization problem (\ref{eq:discreto}).\\
A common approach to solve the generic
constrained optimization problem
\begin{equation}\label{eq:constr}
\begin{array}{l}
    \min_\mathbf{w} J(\mathbf{w}) \\
    \mathrm{s.t.} \\
    H(\mathbf{w})=0,
\end{array}
\end{equation}
\noindent where $J(\mathbf{w})$ is convex and $ H(\mathbf{w})$ is
non-negative, convex and differentiable, is
to convert it into an unconstrained optimization problem.
One way to do this is to use a penalty/continuation method, according to which the constrained
problem \eqref{eq:constr} is approximated by problems of the form:
\[
\min_\mathbf{w} J(\mathbf{w}) + \lambda_k H(\mathbf{w})
\]
where $\{\lambda_k\}$ is an increasing sequence.
In many problems it is necessary
to choose very large values of $\lambda_k$ and it makes the numerical solution process
an extremely difficult one.
Bregman iteration can alternatively be used; it allows one to fix the
value of $\lambda$.\\
One of the central concepts of Bregman iteration is the {\em Bregman distance} \cite{Bregman}
associated with $J$ at point $\mathbf{v}$, defined as:
\begin{equation}\label{eq:distance}
D_J^\mathbf{p}(\mathbf{w},\mathbf{v})=J(\mathbf{w})-J(\mathbf{v})- < \mathbf{p}, \mathbf{w}-\mathbf{v} >,
\end{equation}
\noindent where $\mathbf{p} \in \partial J(\mathbf{v})$
is a subgradient in the subdifferential of $J$ at point $\mathbf{v}$
and $< .,. >$ denotes the canonical vector inner product.\\
By using (\ref{eq:distance}), Bregman iteration can be employed to
reduce (\ref{eq:constr}) in a
short sequence of unconstrained problems,
according to the following iterative scheme:
\begin{equation}\label{eq:bregman1}
\left \{
\begin{array}{l}
   \mathbf{w}_{k+1} = \mathrm{argmin}_\mathbf{w} D_J^\mathbf{p_k}(\mathbf{w}, \mathbf{w}_{k}) +
  \lambda H(\mathbf{w}),\\
   \mathbf{p}_{k+1} =  \mathbf{p}_{k}-\lambda \nabla H(\mathbf{w}_{k+1})
   \in \partial J( \mathbf{w}_{k+1}),\\
\end{array}
\right .
\end{equation}
\noindent with $\lambda>0$.
Inspired by the results obtained in the single case \cite{Noi}, in this section we show that the Bregman iteration is a simple and
very efficient method for solving the portfolio selection problem
(\ref{eq:discreto}), reviewed as a special
case of (\ref{eq:constr}), with $J$ non-smooth convex functional.
In order to derive the expression of $J$ and $H$ for our financial problem,
we introduce two block matrices. Let $N = m \cdot n$ and
$C = diag(C_1, C_2, \ldots, C_m) \in \Re^{N \times N }$
a $m \times m$ diagonal block matrix
with diagonal blocks formed by covariance matrices estimated
at the dates of rebalancing.
$ A \in  R^{(m+1) \times N}$ is the lower bi-diagonal equality constraints
block matrix of dimension $(m+1) \times m$, with diagonal blocks
$A_{i,i}= -\mathbf{1}_n^T$, for $i=1,...,m$, and sub-diagonal blocks
$A_{i+1,i} =( \mathbf{1}_n +\mathbf{r} _{i-1})^T $, for
$i=1,...,m-1$.
The portfolio selection problem (\ref{eq:discreto}) can then be expressed as
(\ref{eq:constr}) with
\begin{equation}
\label{eq:J}
 J(\mathbf{w}) = \mathbf{w}^T C \mathbf{w} + \tau ||\mathbf{w}||_1 ,
\end{equation}
\noindent and
\begin{equation}
\label{eq:H}
 H(\mathbf{w})=\frac{1}{2}||A\mathbf{w}-\mathbf{b}||^2_2,
\end{equation}
\noindent with $\mathbf{b}=(\xi_{init},0,...0,\xi_{term})^T\in \Re^{m+1}$.\\
The choice of the regularization parameter in \eqref{eq:J} plays a key role
to obtain optimal portfolios that meet certain financial requirements and
fidelity to data.
In the single case, where the weights normalization is assumed,
the regularization parameter controls both the sparsity and the
number of short positions;
indeed it can be easily shown that the $l_1$-term permits to penalize the short positions \cite{brodie}:
 \[ \tau ||\mathbf{w}||_1  = 2\tau \sum_{i:w_i<0} |w_i| + \tau.  \]
Unlike the single case, due to self-financial constraint,
the regularization parameter $\tau$ cannot drive
the number of short positions; then it controls only the sparsity. However,
as shown in the following, we observe
that the number of short positions decreases with respect to $\tau$.
For this reason we extend the algorithm proposed in \cite{Noi} to the multi-period case;
the algorithm is based on a modified Bregman iteration to automatically select
the regularization parameter so to satisfy desired financial requirements.
The basic idea is to generate an increasing sequence of parameter values
$\tau_k \in ]0,\tau_{max}]$ that tries
to produce solutions satisfying a fixed financial target, defined
in terms of sparsity or short-controlling.
Given $(\mathbf{w}_{k+1},\mathbf{p}_{k+1})$ provided by \eqref{eq:bregman1} applied to
\[
J_k(\mathbf{w})=\mathbf{w}^T C \mathbf{w} + \tau_k ||\mathbf{w}||_1,
\]
if $\mathbf{w}_{k+1}$ does not satisfy the financial requirement,
$\tau_{k+1} > \tau_k$ is defined.
Consequently, $\mathbf{p}_{k+1}$ must be changed to guarantee the well-definiteness of Bregman iteration.
Following \cite{Noi},
it is possible to show that the vector $\tilde{\mathbf{p}}_{k+1} \in \partial J_{k+1}(\mathbf{w}_{k+1})$,
with
\[ \tilde{\mathbf{p}}_{k+1} = \frac{\tau_{k+1}}{\tau_k}\mathbf{p}_{k+1}+ 2\Bigl(1-\frac{\tau_{k+1}}{\tau_k}\Bigr) C \mathbf{w}_{k+1}
. \]
The proposed multi-period algorithm is summarized in Algorithm~\ref{alg:algo}.
\begin{algorithm}
\caption{Modified Bregman Iteration for portfolio selection}
\label{alg:algo}
\begin{algorithmic}

\State Given $\tau_0>0$, $\tau_{max}$, $\lambda$, $ \theta>1$ \% Model parameters
\State Given $n_{short}$, $n_{act}$ \% Financial target parameters
\State{$k:=0$}
\State{$ \mathbf{w}_0:= \mathbf{0},  \mathbf{p}_0:= \mathbf{0}, \tau_{-1} := \tau_0, $}
\While{``stopping rule not satisfied''}
\State{ $\mathbf{p}_k =  \frac{\tau_{k}}{\tau_{k-1}} \mathbf{p}_k +
\Bigl(1-\frac{\tau_{k}}{\tau_{k-1}}\Bigr) C \mathbf{w}_k$}
\State{ $\mathbf{w}_{k+1} = \mathrm{argmin}_\mathbf{w} D_{J_k}^\mathbf{p_k}(\mathbf{w}, \mathbf{w}_{k}) +
\frac{\lambda}{2} \| A\mathbf{w} - \mathbf{b} \|_2^2$}
\State{ $\mathbf{p}_{k+1} = \mathbf{p}_{k}-\lambda A^T(A  \mathbf{w}_{k+1} -b)$}
\State{$W^-_{k+1} = \{ i   : (\mathbf{w}_{k+1})_i< 0\}$}
\State{$W^a_{k+1} = \{ i   : (\mathbf{w}_{k+1})_i \not = 0\}$}
\If {$|W^-_{k+1}|> n_{short}$ or $|W^a_{k+1}| > n_{act}$}
\State $\eta_{k+1} =\theta$
\Else
\State $\eta_{k+1} = 1$
\EndIf
\State{ $\tau_{k+1} = \min\{\eta_{k+1}\tau_{k},\tau_{max}\}$ }
\State{ $k := k+1$}
\EndWhile
\State{ $\tau_f := \tau_{k+1}$}
\end{algorithmic}
\end{algorithm}
Under suitable hypotheses the convergence of the sequence
$\{ \mathbf{w}_k \}$ to a solution of the constrained problem (\ref{eq:constr})
is guaranteed in a finite number of steps \cite{Osher1}.
Note that the convergence results for Bregman method guarantee the monotonic decrease of
$\| A\mathbf{w}_k- \mathbf{b} \|_2^2$, thus for large $k$ the constraint conditions
are satisfied to an arbitrary high degree of accuracy.
This yields a natural stopping criterion according to a discrepancy principle.\\
Since there is generally no explicit expression for the solution of the sub-minimization
problem involved in (\ref{eq:bregman1}), at each iteration the solution is computed
inexactly using an iterative solver. At this purpose,
we focus on first order methods, which are gradient-based
that converge rather slowly;
however, for large problem dimensions, usually a fast lower-precision solution
is favoured. In particular,
we use the Fast Proximal Gradient method with backtracking stepsize rule (FPG) \cite{fista},
an accelerated variant of Forward Backward algorithm,
suitable for minimizing convex objective functions given by summation of
smooth and non-smooth terms.
Using (\ref{eq:J}), the objective function in the sub-minimization
problem involved in (\ref{eq:bregman1}),
\[\Phi_k(\mathbf{w})=\mathbf{w}^T C \mathbf{w} - <\mathbf{p}_k, \mathbf{w}> + \tau_k ||\mathbf{w}||_1
+ \frac{\lambda}{2} \| A\mathbf{w} - \mathbf{b} \|_2^2, \]
\noindent can be split into the sum of a smooth ($f_k$) and a non-smooth ($g_k$) term
 with
\[ f_k(\mathbf{w})=\mathbf{w}^T C \mathbf{w} - <\mathbf{p}_k, \mathbf{w}> + \frac{\lambda}{2} \| A\mathbf{w} - \mathbf{b} \|_2^2, \;\;\;
 g_k(\mathbf{w})= \tau_k ||\mathbf{w}||_1 . \]
FPG produces a new approximation according to:
\[ \mathbf{w}_{k+1} = prox_{\beta_k g_k} (\mathbf{w}_k-\beta_k \nabla f_k(\mathbf{w}_k)),\]
\noindent with a suitable $\beta_k$ \cite{fista}.
The proximal operator of $g_k$ can be decomposed into a sum of
Soft threshold operator, defined as
\[ Prox_{\tau_k ||\mathbf{w}||_1} (w_i)=sgn(w_i) \left ( |w_i|-min \{|w_i|,\tau_k\} \right ). \]

\section{Numerical results}
\label{sec:4}
In this section, we present the results of some tests
to show the effectiveness of Algorithm \ref{alg:algo} for
solving the multi-period portfolio optimization problem \eqref{eq:discreto}.\\
In Algorithm \ref{alg:algo} we set $\lambda=1,\; \tau_0= 10^{-5}, \; \tau_{max} = 0.5$ and
$\theta = 1.5$.
Iterations are stopped as soon as $\|A\mathbf{w}_k-b\|_2\leq Tol$ with $Tol=10^{-4}$
that, from the financial point of the view, guarantees constraints at a sufficient accuracy.
The maximum number of allowed
Bregman iterations is set to $100$.
Inner iterations are stopped when the relative difference in Euclidean norm between
two successive iterates is less than $Tol_{Inn}=10^{-5}$.\\
The tests have been performed in Matlab R2018a  environment,
on a PC  with Intel Core i5-8250U processors, running Windows 10 Pro 64 bit.\\
We compare our investment strategy with the so-called $1/n$ strategy, where
one invests the same amount of money in all available assets. The portfolio built
following this strategy is referred to as the \emph{naive} portfolio.
It is a common choice to take it as benchmark since
investors often apply this heuristic as well as simple rule
to allocate their wealth across assets,
perceiving it
as a diversification strategy that allows one to reduce risk \cite{naive}.\\
We assume that the investor has one
unit of wealth at the beginning of the planning horizon,
that is, $\xi_{\mathrm{init}}=1$.
In order to compare optimal portfolio with the naive one, we set as
expected final wealth the expected wealth of the naive one, that is,
$\xi_{\mathrm{term}}=\xi_{\mathrm{naive}}$.
The expected wealth of the naive portfolio is obtained by applying recursively the
$1/n$ allocation rule. At each rebalancing date the wealth is evenly invested among the available
securities, that is:
\[
\xi_{\mathrm{naive}}=\frac{1}{n}\left(\ldots\left(\frac{1}{n}\left(\frac{\xi_{\mathrm{init}}}{n}(\mathbf{1}_n+\mathbf{r}_{1})^T\mathbf{1}_{n}\right)
(\mathbf{1}_n+\mathbf{r}_{2})^T\mathbf{1}_{n}\right)\ldots\right)(\mathbf{1}_n+\mathbf{r}_{m})^T\mathbf{1}_{n}
\]
We show results obtained using three real data sets, described below.
The first and the second data sets
come from Fama and French database\footnote{data available at \\
\texttt{http://mba.tuck.dartmouth.edu/pages/faculty/ken.french/data$\_$library.html$\#$BookEquity}}.
The third case study refers to the EURO STOXX50 Index.

\begin{itemize}
\item \textbf{Test 1: FF48}\\
The first database - denoted with FF48 - contains monthly returns
of 48 industry sector portfolios from July $1926$ to December $2015$.
We simulate investment strategies of length $10, 20$ and 30 years,
with annual rebalancing.
\item \textbf{Test 2: FF100} \\
The second database - denoted with FF100 - contains data of $100$
portfolios which are the intersections of $10$ portfolios formed on size
and 10 portfolios formed on the ratio of book equity to market equity.
Also FF100 contains monthly returns from from July $1926$ to December $2015$.
We consider $96$ portfolios of the $100$ available, selected with a
preprocessing procedure which eliminates elements with highest volatilities.
Also in this case we simulate investment strategies of length $10, 20$ and $30$ years,
with annual rebalancing.\\
\item \textbf{Test 3: EURO STOXX50}\\
The third database - denoted with ES50 - contains the daily returns of stocks
included in the EURO STOXX 50 Index
Europe's leading blue-chip index for the Eurozone.
The index covers 50 stocks from 11 Eurozone countries:
Austria, Belgium, Finland, France, Germany, Ireland, Italy, Luxembourg, the Netherlands,
Portugal and Spain. The dataset contains daily returns for each stock in the index from
January $2008$ to December $2013$. 
For this test case we consider both annual and quarterly rebalancing.
\end{itemize}

In table \ref{tab:sparsity} we show results of tests in which the investor requires that
the number of active positions is at most the $30\%$ of the available positions
during the investment period,
that is $n_{act} = 0.3 N$. It is equivalent to require at least the $70\%$ of sparsity in the solution.
In this case there is no limit on the number of short sells
and $n_{short} = N$ in Algorithm \ref{alg:algo}. \\
Optimal and naive portfolios are compared in terms of risk,
given by \eqref{eq:fo}. In particular the ratio between the estimated risk of the $1/n$
strategy and the optimal strategy is reported.
For each test case we report the period of the investment. Note that in all cases the financial
Crisis period is included in the simulation. For optimal portfolios we also report the output value
$\tau_f$ of the regularization parameter, the number of Bregman iterations,
the percentage of short positions and sparsity.
Results show that the required level of sparsity is achieved in all cases, with a modest
number of Bregman iterations.
The $1/n$ strategy is at least $3.25$ times riskier than the the optimal one, for FF100 the ratio is close to $17$.\\
The percentage of shorting varies from the $1\%$ of $N$ for test SP50  to the $9\%$
of $N$ for test FF48 with $10$-years simulation.\\
We note that for the longest investment periods we have the highest values of $\tau_f$;
consequently, also the number of required Bregman iterations and the sparsity level are slightly higher.

\begin{table}[htbp]
\caption{Results for $n_{act} = 0.3 N, \; n_{short} = N$.
The first column contains the label identifying the test case.
The columns $2-7$ contain respectively the investment period, the final value
of the regularization parameter,
the number of Bregman iterations, the percentage of short positions,
the level of sparsity of solutions, the ratio between the risk of the $1/n$ strategy
and the optimal one.}\label{tab:sparsity}
\centerline{
\begin{tabular}{|l|l|c|c|c|c|c|}
  \hline
TEST	& period                & $\tau_f$      & it    & shorts & sparsity & ratio \\ \hline
FF48	&	July 2005-June 2015	&	$8.65E-04$	&	12	&	9\%	&	$72\%$	&	$5.41$	\\
FF48	&	July 1995-June 2015	&	$1.30E-03$	&	14	&	8\%	&	$70\%$	&	$7.72$	\\
FF48	&	July 1985-June 2015	&	$4.38E-03$	&	17	&	6\%	&	$76\%$	&	$8.83$	\\
FF100	&	July 2005-June 2015	&	$5.77E-04$	&	11	&	8\%	&	$80\%$	&	$5.99$	\\
FF100	&	July 1995-June 2015	&	$8.65E-04$	&	12	&	8\%	&	$77\%$	&	$15.35$	\\
FF100	&	July 1985-June 2015	&	$4.38E-03$	&	17	&	4\%	&	$85\%$	&	$16.84$	\\
ES50	&	Jan. 2008-Dec. 2013 (quarterly)	&	$1.71E-04$	&	8	&	1\%	&	$79\%$	&	$3.25$	\\
ES50	&	Jan. 2008-Dec. 2013 (annual)	&	$1.71E-04$	&	8	&	1\%	&	$82\%$	&	$3.45$	\\
 \hline
\end{tabular}
}
\end{table}
As already pointed out, we observe that the number of short positions decreases with
respect to $\tau$. In order to show this behaviour,
in table \ref{tab:taufisso} we report the number of short positions
in optimal portfolios produced with fixed values of the regularization parameter,
obtained by forcing $\theta = 1$ in Algorithm \ref{alg:algo} for FF48 dataset.
Same results are obtained for the other datasets.
As explained in the previous section,
this motivates the introduction of short-selling control in Algorithm \ref{alg:algo}.
Table \ref{tab:short} shows how the proposed adaptive rule for the selection of the regularization
parameter works when short positions are not allowed. This requirement is imposed by setting
$n_{short} = 0, \; n_{act} = N$.
We note that the algorithm is able to find a solution with the desires financial features
in almost all the test. One short position is produced for FF48 in a $30$-years investment.
This result is consistent with the ones reported in table \ref{tab:taufisso},
where one short position is produced for $\tau= 0.1$,
since the value of $\tau_f$ reaches the maximum allowed $\tau_{max}$, and $\tau_{max}$ is $O(10^{-1})$.
We observe that, in general, no-short selling requires quite large values of $\tau$.
This induces high levels of sparsity with a slight increase of risk,
which is anyway lower then the risk of the $1/n$ strategy.
This suggests that the choice of a suitable $\tau_0$ depends on the financial target,
that is, it could be convenient to set a slightly
larger starting value to reduce the iterations if we require no-short controlling.

\begin{table}[htbp]
\caption{Number of short positions produced with fixed values of $\tau$.}\label{tab:taufisso}
\centerline{
\begin{tabular}{|l|c|c|c|}
\hline
FF48 & \multicolumn{3}{c|}{$\tau$} \\
\hline
& $10^{-5}$ & $10^{-3}$ & $10^{-1}$ \\
\hline
July 2005-June 2015	&	202	&	35	&	0	\\
July 1995-June 2015	&	426	&	87	&	0	\\
July 1985-June 2015	&	670	&	225	&	1	\\
\hline
\end{tabular}
}
\end{table}

\begin{table}[htbp]
\caption{Results for $n_{short} = 0, \; n_{act} = N$.
The first column contains the label identifying the test case.
The columns $2-7$ contain respectively the investment period, the final value
of the regularization parameter,
the number of Bregman iterations, the number of short positions,
the level of sparsity of solutions, the ratio between the risk of the $1/n$ strategy
and the optimal one.}\label{tab:short}
\centerline{
\begin{tabular}{|l|l|c|c|c|c|c|}
  \hline
TEST	& period                & $\tau_f$      & it    & nshorts & sparsity & ratio \\ \hline
FF48	&	July 2005-June 2015	&	$6.57E-03$	&	18	&	0	&	88\% & $2.83$	\\
FF48	&	July 1995-June 2015	&	$2.22E-02$	&	22	&	0	&	89\%	& $2.36$ \\
FF48	&	July 1985-June 2015	&	$5.00E-01$	&	47	&	1	&	97\% & $1.46$	\\
FF100	&	July 2005-June 2015	&	$6.57E-03$	&	18	&	0	&	93\% & $1.78$	\\
FF100	&	July 1995-June 2015	&	$3.35E-02$	&	22	&	0	&	96\% & $1.71$	\\
FF100	&	July 1985-June 2015	&	$5.00E-01$	&	47	&	1	&	98\% & $1.22$	\\
ES50	&	Jan. 2008-Dec. 2013 (quarterly)	& $8.65E-04$ &	12	&	0	&	82\%  & $2.82$ \\
ES50	&	Jan. 2008-Dec. 2013 (annual)	& $5.77E-04$ &	11	&	0	&	77\%  &	$2.93$ \\
 \hline
\end{tabular}
}
\end{table}


Our model does not explicitly take into account
transaction costs; however, producing sparse solutions has an impact on transaction costs as well,
allowing one to reduce them. To see this, we count the
transactions.
According to formulation \eqref{eq:discreto} of the financial problem,
if $(\mathbf{w}_{j})_i \neq (\mathbf{w}_{j+1})_i$ we assume that security $i$ has been bought or sold
in the period $[j,j+1[$. Note that this is a pessimistic estimate of transaction costs
because weights could change also for effect of revaluation.
Let us introduce the matrix $G \in \Re^{n \times m}$, with:
\[
\left\{
\begin{array}{ll}
G_{i,j} = 1 & \mbox{ if } (\mathbf{w}_{j})_i \neq (\mathbf{w}_{j+1})_i \\
G_{i,j} = 0 & \mbox{ otherwise} \\
\end{array}
\right.
\]
for $i=1,\ldots,n$ and $j=1,\ldots,m$.
The number of transactions associated with the optimal strategy is then given by:
\[
\mathcal{T}=\sum_{i=1}^{n}\sum_{j=1}^{m}G_{i,j}.
\]
In table \ref{tab:ctcost} we report the number of transactions $\mathcal{T}$ for all the tests.
We denote with $\mathcal{T}_{naive}$ the number of transactions of the $1/n$ strategy.
$\mathcal{T}_{n_act}$ and $\mathcal{T}_{nshort}$ are, respectively,
the number of transactions of the optimal portfolios when either sparsity or
no-shorting are required.
From the table it is evident that the reduction of transaction costs is significant:
the naive portfolio exhibits a number of transactions that is at least three times the
corresponding value of the optimal portfolios. In general costs associated to the
no-short strategy are lower due to the higher sparsity levels observed before. \\

\begin{table}[htbp]
\caption{Number of transactions for all the tests and for the naive and the optimal portfolios
when either sparsity or
no-shorting are required.}\label{tab:ctcost}
\centerline{
\begin{tabular}{|l|l|c|c|c|}
  \hline
TEST	&   period                          & $\mathcal{T}_{naive}$ & $\mathcal{T}_{nact}$  & $\mathcal{T}_{nshort}$ \\ \hline
FF48	&	July 2005-June 2015	            & 480  &	160	&	69			\\
FF48	&	July 1995-June 2015	            & 960  &	352	&	137		 \\
FF48	&	July 1985-June 2015	            & 1440 &	401	&	70			\\
FF100	&	July 2005-June 2015	            & 960  &	242	&	93		\\
FF100	&	July 1995-June 2015	            & 1920 &	584	&	27		\\
FF100	&	July 1985-June 2015	            & 2880 &	553	&	62			\\
ES50	&	Jan. 2008-Dec. 2013 (quarterly)	& 1100 &	279	&  234   \\
ES50	&	Jan. 2008-Dec. 2013 (annual)	& 300  &	76  &  91  \\
 \hline
\end{tabular}
}
\end{table}

\begin{figure}[ht]
  \begin{minipage}[b]{0.5\linewidth}
    \includegraphics[width=1\linewidth]{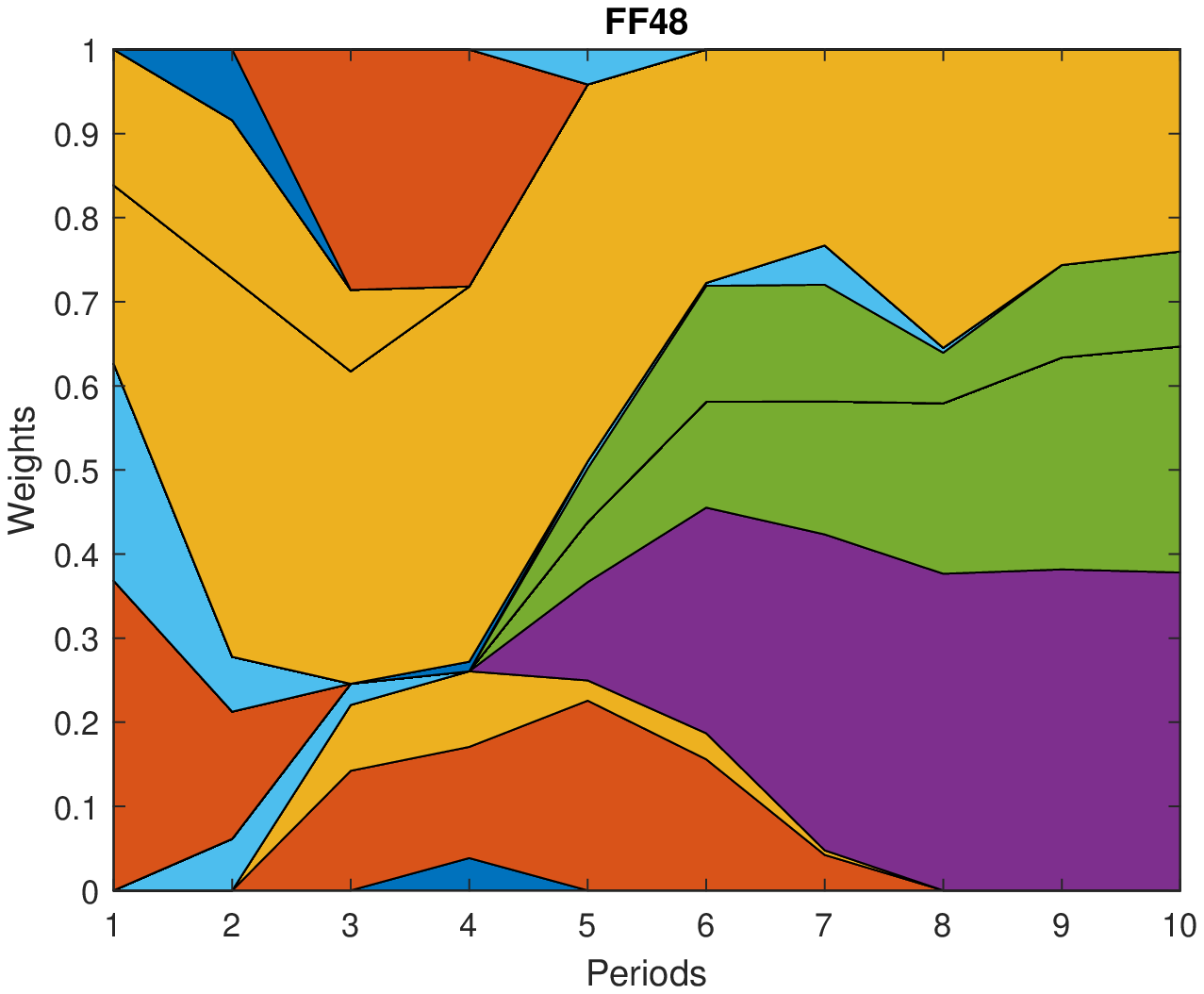}
  \end{minipage}
  \begin{minipage}[b]{0.5\linewidth}
    \includegraphics[width=1\linewidth]{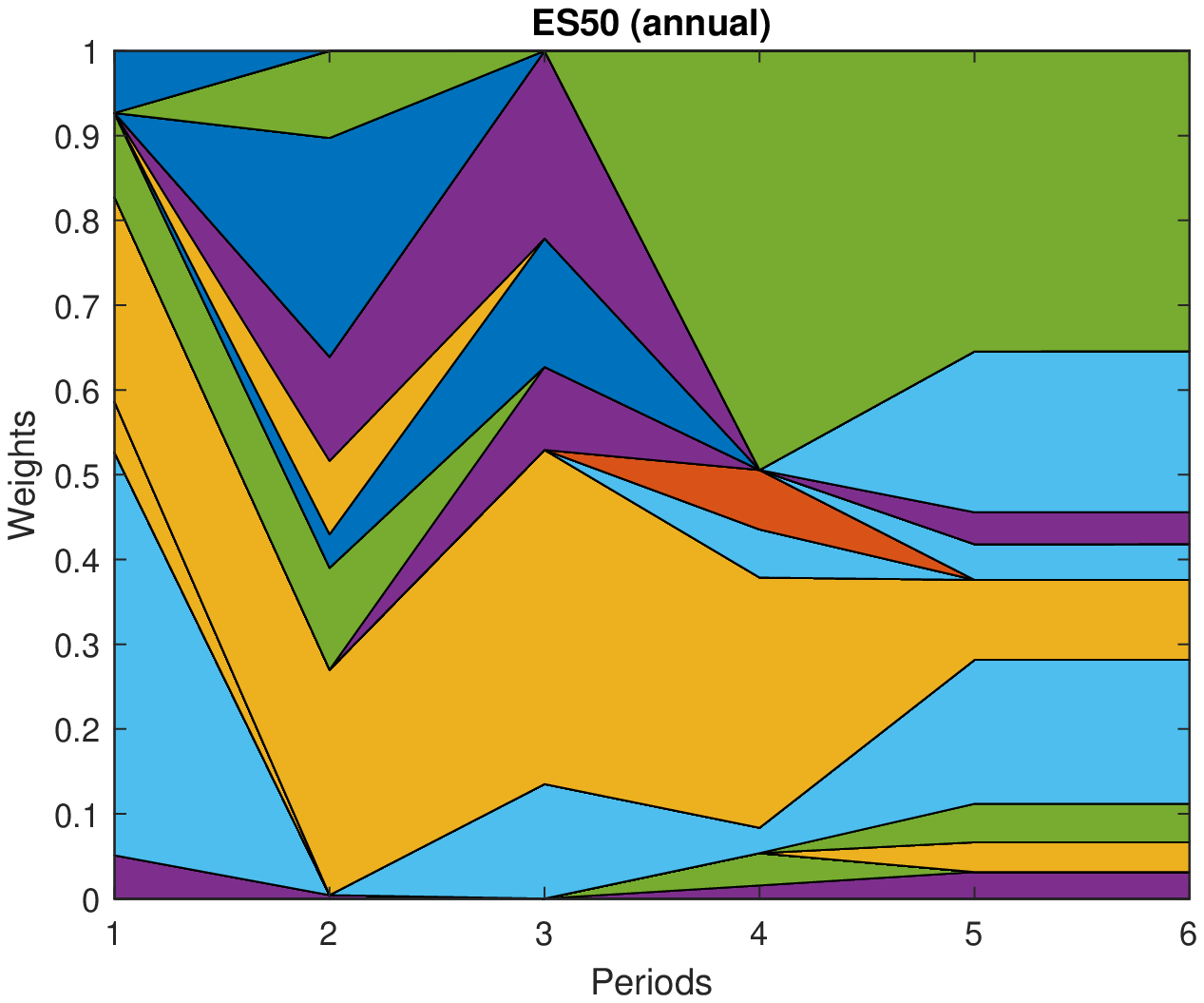}
  \end{minipage}
    \caption{Asset weights over time. Left: FF48, $10$-years
investment; right: ES50 with annual rebalancing. $n_{short}=0$ and $n_{act}=N$
have been set for both.} \label{fig:sparsity}
\end{figure}

Finally, in figure \ref{fig:sparsity} we represent the optimal portfolio weights over time for
two tests; for the sake of readability we consider FF48, $10$-years
investment, and ES50 with annual rebalancing, setting $n_{short}=0$ and $n_{act}=N$ for both.
The number of colored areas at each date represents the number of assets among which the wealth is allocated,
while the height of each colored area at each date represents the amount of wealth allocated in that asset.
We note that only a few assets are involved at each rebalancing date, but the height
of the areas varies over time for most of assets. This suggests that a further
reduction of transaction costs could be obtained, for instance,
modifying the model in such a way to produce sparse solutions
that are kept fixed as much as possible.

\section{Conclusion and future work}
In this work we present a model and the related solution procedure,
based on a modified Bregman iteration,
for the multi-period portfolio selection problem.
The model is corrected with a $l_1$-regularization term to improve conditioning
and obtain sparse solutions.
This has an impact on holding and transaction costs.
A fundamental point is the choice of the regularization parameter that
realizes a good trade-off between sparsity and fidelity to data.
We extend the adaptive rule proposed in \cite{Noi} to the multi-period case.
Numerical results validate our procedure; moreover, we show
that the adaptive selection rule in general can control also the number
of short positions.
Finally we show that transaction costs are reduced, even if they are not
 taken into account in the model. Future work
could concern the introduction of a term that explicitly penalizes transactions
in the model.

%
%

\section{Acknowledgements}
 This work was partially supported by the Research grant of University of Naples ``Parthenope'',
 DR no. 953, November 28th, 2016, and by INdAM-GNCS, under project 2018.

\bibliographystyle{spmpsci}      
\bibliography{biblio}   

\end{document}